\documentclass[12pt]{amsart}
\usepackage{amsfonts,amssymb,amsmath}
\let\germ=\mathfrak
\let\gg=\Gamma
\let\tz=\mathbb

\newtheorem{Theorem}[equation]{Theorem}
\newtheorem{Lemma}[equation]{Lemma}
\newtheorem{Corollary}[equation]{Corollary}
\newtheorem{Proposition}[equation]{Proposition}
\newtheorem{Definition}[equation]{Definition}
\newcommand{\fff}{\ }
\newcommand{\BBB}{B_{SL_2}}
\newcommand{\NNN}{N_{SL_2}}
\newcommand{\ind}{\mathrm{Ind}}
\newcommand{\res}{\mathrm{Res}}

\newcommand{\jad}[2]
{
$\begin{array}{c}
#1\\
#2
\end{array}$
}
\newcommand{\jadv}[2]
{
$\begin{array}{c}
#1\\
\ #2
\end{array}
$
}

\numberwithin{equation}{section}
\begin{document}
\title[Small principal series 
and exceptional duality]{
Small principal series 
and exceptional duality for two simply laced exceptional groups}
\author{Hadi Salmasian}
%\titlerunning{Small principal series and an exceptional duality}
%\institute{
%H. Salmasian \at Queen's University\\
%Department of Mathematics and Statistics\\
%Jeffery Hall, University Avenue\\
%Kingston, Ontario K7L 3N6\\
%Canada\\
%\email{hadi@mast.queensu.ca}
%}
\date{July 19, 2006}
\maketitle
\begin{abstract}
We use the notion of rank defined in 
\cite{salduke}
to introduce and study two correspondences
between small irreducible unitary representations
of the split real simple Lie 
groups of types $\mathbf E_n,\: n\in\{6,7\}$,
and two reductive classical groups.
We show that these correspondences classify all of the 
representations of rank two (in the sense of \cite{salduke}) 
of these exceptional groups.
We study our correspondences for a specific family of degenerate 
principal series representations in detail.
\end{abstract}

%\subjclass[2000]{22E46,22E50,11F27}
%\keywords{Kirillov's orbit method, Mackey analysis, 
%theta correspondence, unitary representations}

%\begin{document}
% \maketitle

\section{Introduction}
\label{introduction}

Construction and classification of small irreducible unitary representations
of non-compact semisimple Lie groups are challenging problems. Small
representations are important because they are natural candidates
for being unipotent representations. In fact many classes of 
small unitary representations are actually automorphic representations.
The most outstanding small unitary representation of semisimple groups 
is probably the oscillator (or the Segal-Shale-Weil) representation.
It happens to be the smallest representation of the metaplectic
group.

Questions about small representations become much
harder for exceptional groups. 
(See \cite{gansavin} and \cite{torasso} where minimal representations
of simple Lie groups are extensively studied.)

This paper is a continuation of the author's work in \cite{salduke}.
In \cite{salduke}, 
the main goal of the author is to define a new
notion of rank for unitary representations of a semisimple 
group over a local field of characteristic zero. 
In principle, this is a generalization of one of 
the main results of 
\cite{licomp} in a fashion that it includes both the classical 
and the exceptional groups at the same time. 
Having an analogous theory of rank, one naturally expects 
that a general classification theorem similar to \cite[Theorem 4.5]{licomp}
should exist. 
Our first main goal is to extend this classification theorem
to exceptional Lie groups. 

Let $G$ be the group of $\tz R$-rational points of a complex, 
absolutely simple,
simply connected algebraic group which is defined over $\tz R$ and is
of type $\mathbf E_n$,
$n\in\{6,7\}$. 
Let $\Pi_2(G)$ denote the set of irreducible unitary representations
of $G$ of rank two (where rank is defined as in 
{\rm \cite[Definition 5.3.3]{salduke}}).  Let $S_1$ be a reductive group given by
\footnote{When $G$ is of type $\mathbf E_6$, the semidirect product $\tz R^\times\ltimes Spin(3,4)$
is explicitly described in the proof of Proposition 
\ref{extensionispossible}.}
$$
S_1=\left\{
\begin{array}{ll}
\tz R^\times\ltimes Spin(3,4)&\textrm{if }G\textrm{ is of type }\mathbf E_6\\
SL_2(\tz R)\times Spin(4,5)&\textrm{if }G\textrm{ is of type }\mathbf E_7
\end{array}
\right.
$$
and let $\Pi(S_1)$ denote the unitary dual of $S_1$.

\begin{Theorem}
\label{main}

There exists an injection 
$\Psi:\Pi_2(G)\mapsto \Pi(S_1)$
which is described in terms of Mackey theory.\vspace{-5mm}\\
\end{Theorem}

The map
$\Psi$, which 
should potentially identify
all of the representations of $G$ of rank two,
is an analogue of the correspondences 
for classical groups which 
appear in \cite{howerank}, \cite{licomp}, and \cite{scaramuzzi}.

From the proof of Theorem \ref{main} it can be seen that one can extend 
the theorem to imply the existence of many maps anlogous to $\Psi$ for several 
other exceptional groups
over complex and $p$-adic fields.
The only difficulties in the proof are certain technicalities with
exceptional groups similar to those dealt with in \cite[\S 5.1,\S 5.2]{salduke}.
However, for the sake of brevity, and because the main point of this paper is
Theorem \ref{main2}, we have stated Theorem \ref{main} for the special case
of real split groups of types $\mathbf E_6$ and $\mathbf E_7$.

Our second main goal is to 
understand the map $\Psi$.
One important common feature of the split forms of types $\mathbf E_6$ and $\mathbf E_7$,
which is one of our motivations for choosing them too, is the existence of certain 
degenerate principal series representations which are of rank two.
See \cite{weissman}, where the author's
motivation for choosing these groups is somewhat similar.
Reducibility points and unitarizability of
this family of principal series representations 
which is introduced in
(\ref{principform}) are studied extensively in \cite{sahiinv} and \cite{barchinietal}.
We use the results in \cite{sahiinv} and \cite{barchinietal} to prove the next two theorems.

\begin{Theorem}
\label{main2}
Let $G$, $\Pi_2(G)$, $S_1$, $\Pi(S_1)$, and $\Psi$ be as above.
Let $P_{\mathrm{ab}}$ be a parabolic subgroup of $G$ whose unipotent radical is
abelian. Let $\chi$ be a unitary character of $P_{\mathrm{ab}}$. If $\pi_\chi$ is the
degenerate principal series representation of $G$ induced from $\chi$,
then $\pi_\chi\in\Pi_2(G)$ and
\begin{itemize}
\item[1.] When $G$ is of type $\mathbf E_6$, $\Psi(\pi_\chi)$ is a unitary character of $S_1$.
\item[2.] When $G$ is of type $\mathbf E_7$, $\Psi(\pi_\chi)$ is isomorphic to
the tensor product of a principal series representation of
$SL_2(\tz R)$ and the trivial representation of $Spin(4,5)$. 
\end{itemize}
\end{Theorem}

Theorem \ref{main2} can also be easily extended to split groups of
types $\mathbf E_6$ and $\mathbf E_7$ over complex and $p$-adic fields.

\begin{Theorem}
\label{main3}
Let $G$ be as above with $n=7$ (i.e., $G$ is of type $\mathbf E_7$).
Let $\Pi_2(G)$ be as above.
Let $I_P(s)$ be defined as in (\ref{principform}).
For any real number $s$ such that $0\leq s< 1$,
Let $\pi_s$ be the complementary series representation
of $G$ 
obtained from $I_P(s)$ as described in 
{\rm\cite[Prop. 7.8, Part (2)]{barchinietal}}. Let
$\pi^\circ$ be 
the representation of $G$ obtained from the irreducible unitarizable
subquotient of $I_P(1)$ as described in {\rm\cite[Prop. 7.8, Part (2)]{barchinietal}}.
Then $\pi^\circ$ and the $\pi_s$ belong to $\Pi_2(G)$.
Moreover, $\Psi(\pi^\circ)$ is the trivial representation of $SL_2(\tz R)\times Spin(4,5)$.
\end{Theorem}

{\noindent\bf Remark.} The reader should note that:

1. When $G$ is of type $\mathbf E_7$, 
it is an interesting problem to understand the ``inverse image'' of 
the unitary dual of $SL_2(\tz R)$ under the map $\Psi$.
For the principal series this follows from Theorem \ref{main2}
and for the trivial representation this follows from Theorem \ref{main3}.
For any $s$ such that $0\leq s<1$, $\Psi(\pi_s)$ should be the tensor
product of a complementary series representation of $G$ and the trivial 
representation of $Spin(4,5)$. For the case of the discrete series of $SL_2(\tz R)$,
Nolan Wallach has suggested a method of construction of
small representations
of $G$ 
based on the method of transfer \cite{wallach}. 
The cases of the complementary series and the discrete series 
will hopefully be addressed 
in a subsequent paper.

2. We could have used
the fancier language of Jordan triple systems to work 
with groups more coherently and probably include some
classical groups as well. (Similar successful attempts 
of using Jordan algebras along these lines 
were made for example in \cite{barchinietal}, \cite{sahiinv},
and \cite{sahidvorsky}.)
However, for classical groups
our results are not new, and for the exceptional groups
that we are going to study it is easier to do things
more explicitly.

This paper is organized as follows. Section 2 merely
introduces our notation and recalls several
facts about representations. Section 3 is devoted to recalling 
the notion of rank in \cite{salduke} and 
the main result of \cite{salduke} on rank.
Section 4 contains a suitable adaptation of a 
result of Scaramuzzi's \cite{scaramuzzi}.
In section 5 we prove Theorem \ref{main} and describe the 
correspondence $\Psi$. In section 6 we prove that the $\pi_\chi$,
the $\pi_s$, and $\pi^\circ$ are of rank two.
(For a description of these representations, see
section 2.)
Sections 7, 8, and 9 are devoted to finding an explicit description
of $\Psi(\pi_\chi)$ and $\Psi(\pi^\circ)$. Section 10 contains
tables concerning roots systems of types $\mathbf E_6$ and $\mathbf E_7$.\vspace{1mm}\\
\noindent{\bf Acknowledgement.}
I would like to thank Roger Howe and Siddhartha Sahi 
for helpful communications, Nolan Wallach for his interest in this
direction of research and illuminating conversations in Toronto, 
and Leticia Barchini and Roger 
Zierau for showing me their joint 
work with Mark Sepanski \cite{barchinietal}. 
I thank the organizing committee 
of the 
conference on representation theory of real reductive groups 
in June 2006 at Snowbird, UT, where some of the proofs
of this paper were simplified.

\section{Notation and preliminaries}
\label{notation}
Let $\tz G$ be a complex, simply connected, absolutely 
simple algebraic group of type $\mathbf E_n$, $n\in\{6,7\}$,  
which is
defined and split over $\tz R$, and let $G$ be the group of
$\tz R$-rational points of $\tz G$.
Let $\germ g$ denote the Lie algebra of $G$. Let
$K$ be the maximal compact subgroup of $G$.

Let $\tz A$ be a maximal split torus of $\tz G$ which is
defined over $\tz R$, and let $A=\tz A\cap G$.
Choose a 
system of roots $\Delta$ associated to $\tz A$ for $\tz G$.
It will induce a root system of $G$.
Choose a positive system $\Delta^+$, and 
let $B$ be the corresponding Borel subgroup of 
$G$. Let $N_B=[B,B]$. 
Let $\Delta_B=\{\alpha_1,...,\alpha_n\}$ 
be a basis for $\Delta^+$. The labelling of the Dynkin diagram
of $\mathbf E_n$ by the $\alpha_i$'s is compatible
with those given in \cite[Planches]{bourbaki} and \cite[Appendix C]{knappbeyond}.
This labelling can be described by the 
diagrams in Figure 1 below. 

\begin{center}
$\begin{array}{c}
6-5-4-3-1\\
|\\
2
\end{array}$
\qquad\qquad
$\begin{array}{c}
{7-6-5-4-3-1\ \ \ \ \ }\\
|\\
2
\end{array}
$

\vspace{-1mm}
{\bf Figure 1.}\vspace{-2mm}
\end{center}

Let $\germ a_\tz C$ be the Lie algebra of $\tz A$.
For any $i\in\{1,...,n\}$, let $\varpi_i\in\mathrm{Hom}_{\tz C}(\germ a_\tz C,\tz C)$ 
denote the fundamental 
weight corresponding to the node labelled by $\alpha_i$, and
let $e^{\varpi_i}$ denote the corresponding character of $\tz A$.
For any $a\in A$, we have $e^{\varpi_i}(a)\in\tz R^\times$. 
Any $a\in A$ is uniquely identified by the values 
of the $e^{\varpi_i}(a)$.

Let $P$ be the standard maximal parabolic
subgroup of $G$ which corresponds to the node
labelled by $\alpha_n$. Let the standard Levi factorisation of $P$ be
$P=L\ltimes N$. The group $N$ is commutative. Let $Q=M\ltimes H$ be the 
standard Levi factorisation of the standard Heisenberg parabolic subgroup of $G$ (where $H$
is a Heisenberg group) and let
$R=S\ltimes U$ be the standard parabolic of $G$ such that
the root system of 
$[S,S]$ is of type $\mathbf D_4$ when $G$ is of type $\mathbf E_6$ and of type
$\mathbf A_1\times \mathbf D_5$ when $G$ is of type $\mathbf E_7$, respectively.
Note that the group $U$ is two-step nilpotent.
Let $P_\gg=L_\gg\ltimes N_\gg$ 
be the standard parabolic introduced in
\cite[\S 3.2]{salduke}. When $G$ is of type $\mathbf E_6$, $[L_\gg,L_\gg]$ is equal to 
the $SL_2(\tz R)$
which corresponds to $\alpha_4$, and when  $G$ is of type $\mathbf E_7$
it is equal to the product of the $SL_2(\tz R)$'s which correspond 
to $\alpha_2,\alpha_3,\alpha_5$ and $\alpha_7$. Let $\beta_1$ be the 
highest root in $\Delta^+$, and let $\beta_2$ be the highest root of
the root system of [$M,M]$. (Obviously we are using the positive 
system for $[M,M]$ which is induced by $\Delta^+$.) 
For every $\alpha\in\Delta$, let $\germ g_\alpha$ be the 
(one-dimensional) root space
of $\germ g$ corresponding to $\alpha$.

We denote the Lie algebras of the groups which appear in the
previous paragraph by 
$\germ p,\germ l,\germ n,\germ q,\germ m,\germ h,\germ r,$
$\germ s,\germ u,\germ p_\gg,\germ l_\gg$ and $\germ n_\gg$
respectively.

Fix a nontrivial positive multiplicative character 
$e^{\Lambda_0}$ of
$L$. For convenience, when $G$ is of type $\mathbf E_7$, we assume $\Lambda_0$ is
the linear functional given in \cite[Definition 2.4]{barchinietal}.
One can extend this character trivially on $N$ to a character
of $P$.
Let $\delta_P$ denote the modular function of $P$.
For any $s\in\tz C$ one can define a degenerate principal series
representation $I_P(s)$ of $G$ as follows:
\begin{equation}
\label{principform}
I_P(s)=\{f\in C^\infty(G)|
\forall p\in P,\forall g\in G: f(gp)=\delta_{P}(p)^{-{1\over 2}}
e^{-s\Lambda_0}(p)
f(g)\}.
\end{equation}
The inner product induced by the norm 
$||f||_K=(\int_K|f(k)|^2dk)^{1\over 2}$
can be used to complete
$I(s)$ and obtain a Hilbert space $\mathcal H_s$.
The action of $G$ on $I_P(s)$ and on $\mathcal H_s$ is 
by left translation:
$$
(g\cdot f)(g_1)=f(g^{-1}g_1).
$$
Obviously, when $s$ is imaginary,
the representation of $G$ on $\mathcal H_s$ is
a principal series representation of $G$
induced from the unitary character of $P$. Let $\overline{P}$
denote the parabolic opposite to $P$. 
There exists an automorphism $\Upsilon$ of $G$
such that $\Upsilon(P)=\overline P$: when $G$ is of type $\mathbf E_7$,
$\Upsilon$ 
is the conjugation by the 
longest element of the Weyl group and when $G$ is of
type $\mathbf E_6$, $\Upsilon$ is
a composition of this conjugation with a diagram automorphism.
From $\Upsilon(P)=\overline P$ it follows that 
one can find a unitary character $\chi$
of $\overline P$ such that
the representation of $G$ on $\mathcal H_s$ 
is isomorphic
to $\ind_{\overline P}^G\chi$. We denote this representation
by $\pi_\chi$.

When $G$ is of type $\mathbf E_7$, reducibility 
points and unitarizability of subquotients of $I_P(s)$ have been 
addressed in \cite{barchinietal},\cite{sahidvorsky},
and \cite{zhang}.
The picture can be concretely described as follows.
For real values of $s$, something similar
to the Wallach set appears (see \cite{barchinietal}). 
More precisely, for
$0\leq s<1$, $I_P(s)$ is irreducible 
and unitarizable and corresponds to a family of
complementary series representations of $G$, which we denote by
$\pi_s$. 
For 
$s=1, 5$, and $9$,
$I_P(s)$ has an irreducible unitarizable subquotient.
The values $s=9$ and $s=5$ correspond to the trivial and the
minimal representations 
of $G$, respectively. The representation obtained
at $s=1$ is denoted by $\pi^\circ$.

From now on, unless stated otherwise, all representations are unitary.
The trivial representation of any group is denoted by ``1''.
The center of a group or a Lie algebra is 
denoted by $\mathcal Z(\cdot)$. For any Hilbert space 
$\mathcal H$, the algebra of bounded operators from $\mathcal H$
to itself is denoted by $\mathrm{End}(\mathcal H)$.

If $G_1$ and $G_2$ are Lie groups where $G_1$ is a 
Lie subgroup of $G_2$, and if for any $i\in\{1,2\}$, $\pi_i$  
is a unitary representation of $G_i$, then $\res_{G_1}^{G_2}\pi_2$
and $\ind_{G_1}^{G_2}\pi_1$ denote restriction and (unitary)
induction.
When there is no ambiguity about $G_2$, we may use ${\pi_2}_{|G_1}$
instead of $\res_{G_1}^{G_2}\pi_2$.
Throughout this paper, we will use two properties of induction 
and restriction which we would like to remind the reader of.
The first property is Mackey's subgroup theorem, as stated in
\cite{mackey}. The second property is the ``projection formula'',
which states that 
$$
\ind_{G_1}^{G_2}(\pi_1\otimes(\res_{G_1}^{G_2}\pi_2))\approx
(\ind_{G_1}^{G_2}\pi_1)\otimes\pi_2.
$$

If $G_1$ is a group acting on a set $X$, then for any $x\in X$
the stabilizer of $x$ inside $G_1$ is denoted by $\mathrm{Stab}_{G_1}(x)$.

Needless to say, let $\tz R^\times$ and
$\tz R^+$ denote (any group naturally isomorphic to) the multiplicative
groups of
nonzero and positive real numbers, respectively.

\section{Representations of small rank}
\label{representations_of_small_rank}

As shown in \cite[Prop. 3.2.6]{salduke}, 
the group $N_\gg$ is a tower of semidirect products
of Heisenberg groups, i.e., \\
\vspace{-5mm}
$$
N_\gg=N_1\ltimes N_2\ltimes N_3
$$
\vspace{-5mm}\ \\
where the $N_i$'s are Heisenberg groups. Note that
$N_3=H$. Let $\germ n_i$ denote the Lie algebra of $N_i$.
For any $i\in\{1,2,3\}$,
let $\rho_i$ be an arbitrary infinite-dimensional irreducible
unitary representation of $N_{4-i}$. By means of the oscillator 
representation, one can extend each $\rho_i$ to a 
representation $\tilde{\rho}_i$ of $N_\gg$. The recipe for extension is given 
in \cite[\S 4.1]{salduke}. We call any representation of $N_\gg$ of the form
$\tilde{\rho}_1$, $\tilde{\rho}_1\otimes\tilde{\rho}_2$ or 
$\tilde{\rho}_1\otimes\tilde{\rho}_2\otimes\tilde{\rho}_3$ a {\it rankable} representation 
of $N_\gg$ of rank one, two or three, respectively (see \cite[Def. 4.1.1]{salduke}). 
The trivial representation
of $N_\gg$ is said to be rankable of rank zero.
One can see that any rankable representation of $N_\gg$ 
is irreducible.
The following theorem is essential to this work.

\begin{Theorem}
\label{mainofduke}
{\rm (\cite[Theorem 5.3.2]{salduke})}. Let $\pi$ be an irreducible unitary representation of $G$.
Then 
the restriction of $\pi$ to $N_\gg$ is 
supported on 
rankable representations of $N_\gg$ of rank $r$, for a fixed
$r\in\{0,1,2,3\}$, which only depends on $\pi$. 
\end{Theorem}

Using Theorem \ref{mainofduke} one can define a notion of rank for
unitary representations of $G$. 
A representation $\pi$ of $G$ is said to be of rank
$r$ if the restriction of $\pi$ to $N_\gg$ is supported on
rankable representations of $N_\gg$ of rank $r$. Theorem
\ref{mainofduke} implies that for irreducible
representations, rank is well-defined.

The only irreducible representation of
$G$ of rank zero is the trivial representation.

\section{Representations of $GL_m(\tz R)$ of rank one}

In this section we prove a result which provides a suitable adaptation of 
\cite[Theorem II.1.1]{scaramuzzi} and will be used in the proof
of Proposition \ref{vonneumann}. The proof of Proposition \ref{scaramodif}
is lengthy but easy, and could be omitted. However, 
to make this manuscript self-contained, we would like
to give its proof in detail.

Fix an integer $m>4$.
Let $GL_m(\tz R)^+$ denote the component group of $GL_m(\tz R)$.
Let $\pi$ be a unitary representation of $GL_m(\tz R)^+$ whose
restriction to $SL_m(\tz R)$ is
of rank one in the sense of {\rm\cite[Definition 5.3.3]{salduke}}.
Let $Q_m$ be the standard Heisenberg parabolic subgroup of $GL_m(\tz R)$
and
$
Q_m^+=Q_m\cap GL_m(\tz R)^+
$.
\begin{Proposition}
\label{scaramodif}
The von Neumann algebra generated by $\pi(Q_m^+)$ is equal to
the von Neumann algebra generated by $\pi(GL_m(\tz R)^+)$.
\end{Proposition}
\begin{proof}
This proposition follows immediately from Lemma
\ref{scaramodif2} below.
\fff\end{proof}

Obviously $GL_m(\tz R)=\{\pm 1\}\ltimes GL_m(\tz R)^+$ for a 
suitable subgroup $\{\pm 1\}$ of the diagonal matrices. (When $m$ is
odd, the semidirect product is actually a direct product.) Observe that
from 
\cite[\S 6]{salduke} it follows that 
a representation of $GL_m(\tz R)$ is of rank one in the sense of 
\cite{scaramuzzi} if and only if the restriction of $\pi$ to 
$SL_m(\tz R)$ is of rank one in the sense of \cite{salduke}.
From now on, a representation of $GL_m(\tz R)^+$ is said to be of rank 
one whenever its restriction to $SL_m(\tz R)$ is of rank one.

\begin{Lemma}
\label{scaramodif1}
Every irreducible
representation $\sigma$ of rank one of $GL_m(\tz R)^+$ 
is equal to the restriction of a representation of $GL_m(\tz R)$. 
\end{Lemma}
\begin{proof}
If $\sigma$ does not extend to a representation of 
$GL_m(\tz R)$, then
$\ind_{GL_m(\tz R)^+}^{GL_m(\tz R)}\sigma$ will be an irreducible 
representation of $GL_m(\tz R)$ of rank one. On the other hand,
the representation \\
\vspace{-3mm}
$$
\res_{GL_m(\tz R)^+}^{GL_m(\tz R)}
\ind_{GL_m(\tz R)^+}^{GL_m(\tz R)}\sigma
$$
\vspace{-5mm}\\
is reducible. (In fact it is a direct sum of $\sigma$ and $\overline\sigma$,
where $\overline\sigma$ is a representation of $GL_m(\tz R)^+$
defined by $\overline\sigma(g)=\sigma(-1\cdot g\cdot -1)$.) But 
from the description of irreducible 
representations of rank one of $GL_m(\tz R)$
in \cite[Theorem II.3.1]{scaramuzzi} as induced representations, 
it follows that the restriction of any irreducible representation of 
rank one of $GL_m(\tz R)$ to $GL_m(\tz R)^+$ remains irreducible. 
Therefore every irreducible representation of $GL_m(\tz R)^+$
of rank one should extend to a representation of $GL_m(\tz R)$, 
which will clearly be of rank one as well.

\fff
\end{proof}

\begin{Lemma}
\label{scaramodif2}
Let $m>4$.
\begin{itemize}
\item[a.] If $\pi$ is an irreducible representation of $GL_m(\tz R)^+$
of rank one, then the restriction of $\pi$ to $Q_m^+$ remains irreducible.
\item[b.] If $\pi,\pi'$ are two distinct irreducible representations of 
$GL_m(\tz R)^+$ of rank one, then their restrictions to $Q_m^+$ are 
nonisomorphic representations.
\end{itemize}
\end{Lemma}
\begin{proof}
Both of these statements follow from \cite[Theorem II.2.1]{scaramuzzi}.
(We advise the reader that in
the notation of \cite{scaramuzzi} 
our $Q_m$ is in fact denoted by $Q_1$. We do not feel obliged to obey 
the notation of \cite{scaramuzzi} since the author does not
use it coherently throughout the paper.)
Let $H_m$ be the unipotent radical of $Q_m$.
Let $J$ be as in \cite[\S I, Equation (18)]{scaramuzzi}
and let $\omega_{tr}$ be the representation of $J$
whose restriction to $H_m$ is the irreducible representation 
$\rho_{tr}$ of \cite[\S II.2]{scaramuzzi} and 
whose extension to $J$ is given 
in \cite[\S II, Equation (39)]{scaramuzzi}.
(We advise the reader that in our situation the indeterminate 
$k$ of \cite[\S II, Equation (39)]{scaramuzzi} is equal to one.)
Let $J^+=J\cap GL_m(\tz R)^+$. Obviously $\omega_{tr}$ can be
considered as a representation of $J^+$ too.
 
If $\pi$ is as in Lemma \ref{scaramodif2}a above, 
then $\pi$ extends to a representation
$\hat{\pi}$
of $GL_m(\tz R)$. By \cite[Theorem II.2.1]{scaramuzzi}, there
exist unitary multiplicative characters $\chi_1,\chi_2$ of $\tz R^\times$
such that \\
\vspace{-3mm}
$$
\res_{Q_m}^{GL_m(\tz R)}\hat{\pi}
=\ind_J^{Q_m}((\chi_1\otimes\chi_2\circ\mathrm{det})\otimes
\omega_{tr}).
$$
(Note that $\chi_1$ and $\chi_2$ play the roles of $\sigma$ and $\chi$
of \cite[Theorem II.2.1]{scaramuzzi}.) 
By Mackey's subgroup theorem one can see that \\
\vspace{-3mm}
$$
\res_{Q_m^+}^{GL_m(\tz R)^+}{\pi}
=\ind_{J^+}^{Q_m^+}((\chi_1\otimes\chi_2\circ\mathrm{det})\otimes
\omega_{tr}).
$$
Recall that $H_m$ is the unipotent radical of $Q_m$. 
The restriction of the representation $\omega_{tr}$ to 
$H_m$ is irreducible. Therefore $\omega_{tr}$ is an irreducible 
representation of $J^+$ too.
Standard Mackey theory tells us that for any irreducible representation
$\nu$ of $J^+/H_m$, $\ind_{J^+}^{Q_m^+}(\nu\otimes
\omega_{tr})$ is an irreducible representation of $Q_m^+$. Therefore
$\pi_{|Q_m^+}$ is irreducible. This proves Lemma \ref{scaramodif2}a.

Next we prove Lemma \ref{scaramodif2}b. By Lemma
\ref{scaramodif1} 
it suffices to show that if $\hat{\pi}$ and 
$\hat{\pi}'$ are two
rank one representations of $GL_m(\tz R)$ whose restrictions to 
$Q_m^+$ are isomorphic, then the restrictions of $\hat{\pi}$ 
and $\hat{\pi}'$
to $GL_m(\tz R)^+$ are isomorphic as well. Suppose 
$\chi_1,\chi_2,\chi_1',\chi_2'$ are unitary multiplicative characters
of $\tz R^\times$ such that\\
\vspace{-2mm}
$$
\res_{Q_m}^{GL_m(\tz R)}\hat{\pi}
=\ind_J^{Q_m}((\chi_1\otimes\chi_2\circ\mathrm{det})\otimes
\omega_{tr})
$$\vspace{-5mm}\\
and\vspace{-2mm}
$$
\res_{Q_m}^{GL_m(\tz R)}\hat{\pi}'
=\ind_J^{Q_m}((\chi_1'\otimes\chi_2'\circ\mathrm{det})\otimes
\omega_{tr}).
$$
Obviously we have
\vspace{-3mm}
$$
\res_{Q_m^+}^{GL_m(\tz R)}\hat{\pi}
=\ind_{J^+}^{Q_m^+}((\chi_1\otimes\chi_2\circ\mathrm{det})\otimes
\omega_{tr})
$$
\vspace{-7mm}\\
and \vspace{-1mm}
$$
\res_{Q_m^+}^{GL_m(\tz R)}\hat{\pi}'
=\ind_{J^+}^{Q_m^+}((\chi_1'\otimes\chi_2'\circ\mathrm{det})\otimes
\omega_{tr}).
$$
By standard Mackey theory, if
the restrictions of $\hat{\pi}$ and $\hat{\pi}'$ to
$Q_m^+$ are isomorphic, then we have\\
\vspace{-5mm}
\begin{equation}
\label{characterequal}
\chi_1(a)=\chi_1'(a) 
\textrm{ for every } a\in \tz R^\times 
\textrm{ and } \chi_2(a)=\chi_2'(a) 
\textrm{ for every } a\in\tz R^+. 
\end{equation}

Let $P_m$ be the standard parabolic subgroup of $GL_m(\tz R)$
whose Levi factor is 
$$
GL_1(\tz R)\times GL_{m-1}(\tz R)
$$
and
let $P_m^+=P_m\cap GL_m(\tz R)^+$.
If $\chi^1$ and $\chi^2$ are characters of 
$GL_1(\tz R)$ and $GL_{m-1}(\tz R)$ 
respectively, then we define the representation $\chi^1\times\chi^2$
of $GL_m(\tz R)$ induced from $P_m$ 
as in \cite[\S II.3]{scaramuzzi}.
From \cite[Theorem II.3.1]{scaramuzzi} 
it follows that
$$
\hat{\pi}=\ind_{P_m}^{GL_m(\tz R)}(\chi_2^{-1}\chi_1\times \chi_2\circ\det)
\qquad\mathrm{and}\qquad
\hat{\pi}'=\ind_{P_m}^{GL_m(\tz R)}({\chi_2'}^{-1}\chi_1'\times \chi_2'
\circ\det).
$$
Consequently, from Mackey's subgroup theorem it follows that\\
\vspace{-3mm}
$$
\hat{\pi}_{|GL_m(\tz R)^+}=
\ind_{P_m^+}^{GL_m(\tz R)^+}(\chi_2^{-1}\chi_1\times \chi_2\circ\det)
$$ 
\vspace{-5mm}and\\
\vspace{-1mm}
$$
\hat{\pi}'_{|GL_m(\tz R)^+}=
\ind_{P_m^+}^{GL_m(\tz R)^+}({\chi_2'}^{-1}\chi_1'\times \chi_2'\circ\det).
$$
\vspace{-3mm}\\
Equalities in (\ref{characterequal}) imply that
$\chi_2^{-1}\chi_1\times \chi_2\circ\det$ 
and ${\chi_2'}^{-1}\chi_1'\times \chi_2'\circ\det$
are identical characters of $P_m^+$, i.e., restrictions of $\hat{\pi}$ and
$\hat{\pi}'$ to $GL_m(\tz R)^+$ are isomorphic.
This proves Lemma \ref{scaramodif2}b.

\fff
\end{proof}

\section{Proof of Theorem \ref{main}}

It was shown in \cite[Prop. 4]{salmanus} 
that the only irreducible
representation of $G$ of rank one is the {\it minimal} representation
of $G$. Our concentration throughout the rest of this manuscript 
will be on representations of $G$ of rank two.
It was shown in \cite{kazhdansavin} that the minimal representation 
of $G$ is
irreducible when restricted to the Heisenberg parabolic.
Our next task is to prove a similar, but much stronger
version of this fact for irreducible representations of rank two.
Namely, we will show that for irreducible representations of
rank two, the restriction of $\pi$ to the parabolic subgroup $R$
determines $\pi$ uniquely. Our method of proof is an adaptation of an 
idea originally due to Howe \cite{howerank}.

The parabolic $Q$ can be expressed as 
\begin{equation}
\label{qlevifac}
Q=(\tz R^+\times(\{\pm 1\}\ltimes [M,M]))\ltimes H
\end{equation}
where $\tz R^+$ and $\{\pm 1\}$ are appropriate 
subgroups of $A$. (Note the positions of direct and semidirect products.)
The element $-1\in\{\pm 1\}$ corresponds to the element $a\in A$ such that for $G$ is of type $\mathbf E_6$
we have $e^{\varpi_2}(a)=-1$ and $e^{\varpi_j}(a)=1$ for any $j\neq 2$, and
for $G$ is of type $\mathbf E_7$ we have $e^{\varpi_1}(a)=-1$ and $e^{\varpi_j}(a)=1$ for any $j\neq 1$.

Let $\rho$ be any infinite-dimensional irreducible 
unitary representation of $H$.
Note that $\rho$ extends to a representation 
$\hat{\rho}$ of $[Q,Q]$, and in fact this extension is 
unique since $[M,M]$ is a perfect group. 
Since 
the action of $Q$ on the center of $H$ has only one open orbit
and the stabilizer of every point of this orbit is $[Q,Q]$,
it follows 
from \cite[\S 5.2]{salduke} that
one can express the restriction of $\pi$ to $Q$ as \\
\vspace{-4mm}
$$
\pi_{|Q}=\ind_{[Q,Q]}^Q(\nu\otimes \hat{\rho})
$$\vspace{-4mm}\\
where $\nu$ is a representation of $[M,M]=[Q,Q]/H$.
Let $Q^+$ be the component group of $Q$. Then we have \\
\vspace{-4mm}
$$
Q^+=(\tz R^+\times [M,M])\ltimes H.
$$ 
\vspace{-4mm}\\
One can extend $\nu$ to the representation $1\otimes\nu$ 
of $\tz R^+\times [M,M]$ and therefore obtain an extention
$\hat{\nu}$ of $\nu$ to $Q^+$.
(Note that it may not necessarily be possible to
extend $\nu$ to all of $Q$.)
Using the projection formula we have\\\vspace{-3mm}
$$
\ind_{[Q,Q]}^Q(\nu\otimes \hat{\rho})=\ind_{Q^+}^Q
\ind_{[Q,Q]}^{Q^+}(\nu\otimes \hat{\rho})=\ind_{Q^+}^Q
(\hat{\nu}\otimes\ind_{[Q,Q]}^{Q^+}\hat{\rho}).
$$\vspace{-3mm}\\
Let $\tau=\hat{\nu}\otimes\ind_{[Q,Q]}^{Q^+}\hat{\rho}$ and
$\eta=\ind_{[Q,Q]}^{Q^+}\hat{\rho}$.
If $\mathcal H_1$ and $\mathcal H_2$ denote the Hilbert
spaces of the representations $\hat{\nu}$ and $\eta$, 
then $\tau$ is a representation with Hilbert space 
$\mathcal H_1\otimes \mathcal H_2$.

\begin{Lemma}
\label{vonneu1}
The von Neumann algebra generated
by $\tau(Q^+\cap R)$ is equal
to the von Neumann algebra
generated by $\tau(Q^+)$.
\end{Lemma}

\begin{proof}
The key point is that
the restriction of $\eta$ to 
the subgroup $\tz R^+\ltimes H$ of $Q^+$ is irreducible
(see \cite{kazhdansavin}). The idea 
is that the $H$-spectrum of $\eta$ is multiplicity-free and $\tz R^+$ acts
on it transitively.
Since the restriction of $\hat{\nu}$ to $\tz R^+\ltimes H$
is trivial, it follows that the von Neumann algebra generated
by $\tau(\tz R^+\ltimes H)$ contains every operator of the
form $I\otimes T$ inside $\mathrm{End}(\mathcal H_1\otimes\mathcal H_2)$.
Consequently, the von Neumann algebra $\mathcal A$ generated by
$\tau(Q^+\cap R)$ contains all the operators of the form
$\hat{\nu}(q)\otimes T$ for $q\in Q^+\cap R$, where
$T$ is any arbitrary element of 
$\mathrm{End}(\mathcal H_2)$.

Next we observe that $\hat{\nu}$
is actually a representation of rank one of 
the reductive group $M^+=Q^+/H$. 
But for classical groups it 
can be shown that the two notions of rank in \cite{licomp} and
in \cite{salduke} are essentially 
equivalent (see \cite[\S 6]{salduke}). Therefore from 
\cite[Theorem 4.5]{licomp} and Lemma \ref{scaramodif}
it follows that the
von Neumann algebra generated by $\hat{\nu}(Q^+\cap R)$ is equal to
the von Neumann algebra generated by $\hat{\nu}(Q^+)$. 
(When $G$ is of type $\mathbf E_6$, the latter statement follows immediately from Proposition 
\ref{scaramodif}.
When $G$ is of type $\mathbf E_7$, it follows from the fact that the von Neumann algebra
generated by $\hat{\nu}(Q^+\cap R)$ contains the von Neumann algebra 
generated by $$\hat{\nu}([Q,Q]\cap R),$$ which 
by \cite[Theorem 4.5]{licomp} is equal to the von Neumann algebra generated
by $\hat{\nu}([Q,Q])$.)
This implies that
$\mathcal A$ contains all the operators of the form 
$\hat{\nu}(q)\otimes T$ for every $q\in Q^+$, which
proves Lemma \ref{vonneu1}.

\fff\end{proof}

\begin{Lemma}
\label{vonneu22}
The von Neumann 
algebra generated by $\pi(Q)$ is equal to the von Neumann algebra
generated by $\pi(Q\cap R)$.
\end{Lemma}
\begin{proof}
Let $\sigma=\ind_{Q^+}^Q\tau$. We have $\pi_{|Q}=\sigma$. By the double commutant theorem, it
suffices to show that every $Q\cap R$-intertwining operator for $\sigma$
is actually a $Q$-intertwining operator. Let 
$\mathcal H_\tau$ be the Hilbert space of
$\tau$. Recall that $-1\in\{\pm 1\}\subset Q$ (see equation 
(\ref{qlevifac})).
One can realize $\sigma$ on $\mathcal H_\tau\oplus\mathcal H_\tau$
as
\begin{eqnarray*}
\sigma(q)(v\oplus w)&=&\tau(q)v\oplus\overline\tau(q)w\qquad 
\textrm{for any }\ q\in Q^+\\
\sigma(-1)(v\oplus w)&=&w\oplus v\\
\end{eqnarray*}
\vspace{-12mm}\\
where $\overline\tau$ is the representation of $Q^+$ on $\mathcal H_\tau$
obtained by twisting by $-1$, i.e., \\
\vspace{-3mm}
$$
\overline\tau(q)=\tau(-1\cdot q \cdot -1)
\quad\textrm{for any }\ q\in Q^+.
$$ 
\vspace{-3mm}\\
Any element of $\mathrm{End}
(\mathcal H_\tau\oplus\mathcal H_\tau)$ can be represented 
by a matrix 
$$
T=\left[\begin{array}{cc} 
T_1&T_2\\
T_3&T_4
\end{array}
\right]
$$
where $T_1,T_2,T_3,T_4$ are elements of $\mathrm{End}(\mathcal H_\tau)$.
Let $T$ be a $Q\cap R$-intertwining operator.
The fact that $T$ commutes with $\sigma(-1)$ implies that
$T_1=T_4$ and $T_2=T_3$. The fact that $T$ commutes with the action of
$\sigma(q)$ for $q\in Q^+\cap R$ 
implies that $T_2$ is a $Q^+\cap R$-interwining operator
between $\tau$ and $\overline\tau$
and $T_1$ is a $Q^+\cap R$-intertwining operator of 
$\tau$. However, the restrictions of
$\tau$ and $\overline\tau$ to $H$ are disjoint (one of $\tau,\overline\tau$
is supported on representations of $H$ with ``positive'' central character
whereas the other one is supported on representations with ``negative'' 
central
character). Therefore $T_2$ should be zero. On the other hand 
since by Lemma \ref{vonneu1} the von Neumann algebras generated by $\tau(Q^+\cap R)$ and $\tau(Q^+)$ 
are the same, it follows that $T_1$ is a $Q^+$-intertwining operator 
of $\tau$ as well. Consequently, $T$ is of the form
$$
T=\left[
\begin{array}{cc}
T_1&0\\
0&T_1
\end{array}
\right]
$$
which implies that $T$ is a $Q$-intertwining operator for $\sigma$.

\fff
\end{proof}

\begin{Proposition}
\label{vonneumann}

Let $\pi$ be any unitary representation of $G$ of rank two on a Hilbert
space $\mathcal H$.
Then the von Neumann algebra (inside $\mathrm{End}(\mathcal H)$) generated
by $\pi(R)$ is identical to the von Neumann algebra
generated by $\pi(G)$.
\end{Proposition}
\begin{proof}

The von Neumann algebra generated by $\pi(R)$ contains 
the von Neumann algebra generated by $\pi(Q\cap R)$, which by Lemma 
\ref{vonneu22} is equal to 
the von Neumann algebra generated by $\pi(Q)$. Since $Q$ and $R$
are both maximal parabolics, the group generated by them is equal to
$G$. Therefore the von Neumann algebra generated by $\pi(R)$ 
contains the von Neumann
algebra generated by $\pi(G)$.

\fff
\end{proof}
\begin{Corollary}
\label{corondisting}
If $\pi$ is an irreducible representation of $G$ of rank two, 
then the restriction of $\pi$ to $R$ is irreducible 
and uniquely determines $\pi$.
\end{Corollary}

We are now able to apply the standard machinery of 
Mackey \cite{mackey} to our situation. Let $\pi$ be an irreducible 
unitary representation of
$G$ of rank two. 
Recall that the Levi factorisation of $R$ is $R=S\ltimes U$.
The group $[S,S]$ is the 
set of $\tz R$-points of a complex simple simply connected algebraic group 
which is defined and split over $\tz R$ and whose 
root system is of type $\mathbf D_4$ when $G$ is of type $\mathbf E_6$ and of type 
$\mathbf A_1\times \mathbf D_5$ when $G$ is of type $\mathbf E_7$. Therefore 
we have 
$$
[S,S]=\left\{
\begin{array}{ll}
Spin(4,4)&\qquad\textrm{if } G \textrm{ is of type }\mathbf E_6\\
SL_2(\tz R)\times Spin(5,5)&\qquad\textrm{if } G \textrm{ is of type }\mathbf E_7.
\end{array}
\right.
$$
Our next task is to understand the restriction of 
$\pi$ to $U$. 
Recall that $U$ is two-step nilpotent.
Let $\mathcal Z(U)$ denote the center of $U$. 
Let $\sigma$ be an irreducible unitary representation of 
$U$. For every element $z\in\mathcal Z(U)$, $\sigma(z)$ is a
scalar.
If $\mathcal Z(U)\subseteq \ker(\sigma)$ then 
$\sigma$ should be one-dimensional. Now suppose 
$\mathcal Z(U)\nsubseteq\ker(\sigma)$.
The group $\mathcal Z(U)$ is invariant under the action of $S$ and 
one can see that the action of the spin factor of $[S,S]$
on the Lie algebra $\mathcal Z(\germ u)$ is in fact identical to the 
standard representation $\tz R^{n-2,n-2}$ of $Spin(n-2,n-2)$.
When $G$ is of type $\mathbf E_7$, the factor $SL_2(\tz R)$ of $[S,S]$ acts 
on $\mathcal Z(\germ u)$
trivially.

From the existence of a non-degenerate 
$Spin(n-2,n-2)$-invariant bilinear form 
$<\cdot,\cdot>$
on $\tz R^{n-2,n-2}$,
it follows that for any unitary character 
$\chi$ of $\mathcal Z(U)$ there exists an element
$v\in \mathcal Z(\germ u)$
such that  
\vspace{-2mm}
\begin{equation}
\label{characterofUU}
\chi(x)=e^{<v,\log x>\sqrt{-1}}\qquad\textrm{for any }x\in \mathcal Z(U).
\end{equation}
\vspace{-6mm}\\
Here ``$\log x$'' means the inverse of the exponential map
$\exp:\mathcal Z(\germ u)\mapsto \mathcal Z(U)$.
From $\mathcal Z(U)\nsubseteq \ker(\sigma)$ it follows that
$v\neq 0$.
The action of $S$ on $\mathcal Z(\germ u)$ 
has three
orbits: an open orbit, the set of nonzero isotropic vectors,
and the origin (see \cite{howetan}
and \cite{kaneyuki}).
\begin{Definition}
Consider an irreducible unitary representation of $U$ whose
restriction to $\mathcal Z(U)$ acts by a character
of the form given in (\ref{characterofUU}) for some $v\neq 0$. 
We call this representation small if $v$ is a nonzero isotropic vector,
and we call it big if $v$ belongs to the open orbit.
\end{Definition}

Suppose $G$ is of type $\mathbf E_n$. If $\sigma$ is a big representation of $U$ then 
$U/\ker(\sigma)$ is a Heisenberg group of dimension $2^{n-2}+1$,
whereas if $\sigma$ is small then $U/\ker(\sigma)$ is
a direct product of the additive group of 
$\tz R^{2^{n-3}}$ and a Heisenberg
group of dimension
$2^{n-3}+1$.

Recall 
that $\germ u$ and $\germ n_3$ represent the Lie algebras of $U$ and
$N_3$ respectively. As a vector space, we can write
$\germ u$ as a direct sum  
\vspace{-2mm}
\begin{equation}
\label{descripuxyz}
\germ u=\germ X\oplus \germ Y\oplus \mathcal Z(\germ u)
\end{equation}
\vspace{-6mm}\\
where:
\begin{itemize}
\item[a.] Each of $\germ X$ and $\germ Y$ 
is a direct sum of certain root spaces $\germ g_\alpha$.

\item[b.]$\germ X=\germ u\cap \germ l$.

\item[c.] $\mathcal Z(\germ u)$ is the center of $\germ u$.

\end{itemize}
Note that these conditions identify $\germ Y$ uniquely. For an explicit description
of this decomposition, see section \ref{Appendix}.
Similarly, we can write $\germ n_3$ as a direct sum
\vspace{-2mm}
$$\germ n_3=\germ W\oplus\germ W^*\oplus \mathcal Z(\germ n_3)$$
\vspace{-6mm}\\
where:

\begin{itemize}
\item[a.] Each of $\germ W$ and $\germ W^*$ is a 
direct sum of certain root spaces $\germ g_\alpha$.

\item[b.] $\germ W=\germ n_3\cap \germ l$.

\item[c.] $\mathcal Z(\germ n_3)$ is the center of $\germ n_3$.

\end{itemize}
Again $\germ W^*$ is uniquely identified (see section \ref{Appendix}). In fact $\germ W$ 
and $\germ W^*$ correspond to a polarization of the symplectric 
vector space $\germ n_3/\mathcal Z(\germ n_3)$.

\begin{Lemma}
\label{basislemma}
The Lie algebra $\mathcal Z(\germ u)$ contains 
$\germ g_{\beta_1}$ and $\germ g_{\beta_2}$.
For any root $\alpha\notin\{\beta_1,\beta_2\}$, 
if $\germ g_\alpha\subseteq \mathcal Z(\germ u)$
then $\germ g_\alpha\subseteq\germ W^*$. Moreover,
there exist bases 
$$
\{e_1,...,e_{n-2},e_{-1},...,e_{-(n-2)}\}\textrm{ and  }
\{f_2,...,f_{n-2},f_{-2},...,f_{-(n-2)}\}
$$ 
of $\mathcal Z(\germ u)$ and 
$\germ W\cap \germ s$ respectively,
such that

\begin{itemize}
\item[a.] For any $i\in\{-(n-2),...,n-2\}$, 
there exist $\alpha_1,\alpha_2\in\Delta^+$ such that 
$e_i\in\germ g_{\alpha_1}$ and $f_i\in\germ g_{\alpha_2}$.

\item[b.] $e_1\in\germ g_{\beta_1}$ and $e_{-1}\in\germ g_{\beta_2}$.

\item[c.] For any $i,j\in\{2,...,n-2\}$,
$[e_{\pm i},f_{\pm j}]=\delta_{\pm i,\pm j}e_1$.

\item[d.] For any $i,j\in\{1,...,n-2\}$, we have 
$$
<\!\!e_i,e_j\!\!>=<\!\!e_{-i},e_{-j}\!\!>=0
\qquad\textrm{and}\qquad <\!\!e_i,e_{-j}\!\!>=\delta_{ij}.
$$

\item[e.] For any root $\alpha$, if $\germ g_\alpha\subseteq \germ n_3$
and $\germ g_\alpha\nsubseteq\mathrm{Span}_{\tz R}\{e_i,f_i\ |\ -(n-2)\leq i 
\leq n-2 \}$ then $[e_{-1},\germ g_\alpha]=0$.

\end{itemize}
\end{Lemma}
\begin{proof}
See section \ref{Appendix}.

\fff
\end{proof}

\begin{Proposition}
\label{ranktwoopenorbit}
Let $\pi$ be an irreducible representation of $G$ of rank two.
Then the restriction of $\pi$ to $U$ is supported on
big representations of $U$. 
 
\end{Proposition}
\begin{proof}
Since $\mathcal Z(U)\subset N_\gg$, 
it suffices to prove that the restriction of
any rankable representation of $N_\gg$ of rank two
to $\mathcal Z(U)$ is a direct integral of characters
of the form given in equation (\ref{characterofUU}) for which
$v$ belongs to the open $S$-orbit of $\mathcal Z(\germ u)$.
Let $\tilde{\rho}_1$ be a rankable representation of $N_\gg$ of
rank one, obtained by extending a representation of $N_3=H$.
By Lemma \ref{basislemma} for any $i,j\in\{\pm 2,\ldots,\pm (n-2)\}$
we have 
$$<\!e_j,[e_{-1},f_i]\!>=-<\!e_j,[f_i,e_{-1}]\!>=<\![f_i,e_j],e_{-1}\!>=
<\!-\delta_{i,j}e_1,e_{-1}\!>=-\delta_{i,j}$$
which implies that $[e_{-1},f_i]=-e_{-i}$. 
Using 
Lemma \ref{basislemma}e
and the formulas which describe the oscillator 
representation (see \cite[\S 1]{howerank}), 
one can see that the
restriction of $\tilde{\rho}_1$ to $\mathcal Z(U)$ is 
a direct integral of characters of the form
given in (\ref{characterofUU}) for 
\vspace{-3mm}
$$
v=te_1+\sum_{i=2}^{n-2} (ta_ie_i+ta_{-i}e_{-i})-
(\sum_{i=2}^{n-2}a_ia_{-i}t)e_{-1}
$$
\vspace{-4mm}\\
where $t$ and the $a_i$'s are real numbers and $t\neq 0$.
Obviously $<v,v>=0$, i.e., $v$ is an isotropic vector.
From Lemma \ref{basislemma}e and the oscillator
representation it also
follows that the restriction of a rankable representation of
rank two of $N_\gg$ to $\mathcal Z(U)$ is a direct integral of
characters for which $v$ is given by
\vspace{-2mm}
$$
v=te_1+\sum_{i=2}^{n-2} (ta_ie_i+ta_{-i}e_{-i})-
(s+\sum_{i=2}^{n-2}a_ia_{-i}t)e_{-1}
$$\vspace{-4mm}\\
where $t$ and $a_i$'s are as above and $s$ is a nonzero 
real number. Obviously $<v,v>\neq 0$.

\fff

\end{proof}

Let $\pi$ be an irreducible representation of $G$ of rank two.
Let $\sigma$ be the irreducible unitary representation of 
$U$ with central character given by (\ref{characterofUU})
where $v=e_1+e_{-1}$.
By Mackey theory, one can write the restriction of $\pi$ to $R$
as
$$
\pi_{|R}=\ind_{R_1}^R\eta
$$
where $R_1=\mathrm{Stab}_S(e_1+e_{-1})\ltimes U$
and $\eta$ is an irreducible representation of $R_1$ 
with the property that for some $m\in \{0,1,2,\ldots,\infty\}$, 
$\res_{U}^{R_1}\eta=m\sigma$.
If $\sigma$ can be extended
to a representation $\hat{\sigma}$ 
of $R_1$,
then $\eta$ can be written as a tensor product
$\eta=\tau\otimes\hat{\sigma}$ where $\hat{\sigma}$ is the
extension of $\sigma$ to $R_1$ and $\tau$ is an irreducible 
representation
of $\mathrm{Stab}_S(e_1+e_{-1})$ which is extended 
(trivially on $U$) to $R_1$. 
The representation $\tau$ is uniquely
determined by $\pi_{|R}$, and 
by Corollary \ref{corondisting}
we obtain an injection 
\begin{equation}
\label{correspondence}
\Psi:\pi\mapsto\Psi(\pi),
\end{equation}
where $\Psi(\pi)=\tau$,
from $\Pi_2(G)$ 
into the unitary dual of $\mathrm{Stab}_S(e_1+e_{-1})$.

The next proposition justifies the existence of the
extension $\hat{\sigma}$.

\begin{Proposition}
\label{extensionispossible}
Let $\sigma$ be the big representation of $U$ associated to
the character given by (\ref{characterofUU}) where
$v=e_1+e_{-1}$. Then $\sigma$ can be extended to a representation
$\hat{\sigma}$ of $R_1$. The extension is unique
when $G$ is of type $\mathbf E_7$.
\end{Proposition}

\begin{proof}
Let $S_1=\mathrm{Stab}_S(e_1+e_{-1})$.
It is well known that for the action of 
$Spin(2k,\tz C)$ on $\tz C^{2k}$, 
if $w\in\tz C^{2k}$ lies outside 
the variety of isotropic vectors, then
\begin{equation}
\label{localstab}
\mathrm{Stab}_{Spin(2k,\tz C)}(w)=Spin(2k-1,\tz C).
\end{equation}
From (\ref{localstab}) and some elementary calculations
it follows that when $G$ 
is of type $\mathbf E_6$, 
$$
S_1=\tz R^\times\ltimes Spin(3,4)
$$
where the element $-1\in\tz R^\times$ is an element of $A$
given by $e^{\varpi_1}(a)=e^{\varpi_6}(a)=-1$ and
$e^{\varpi_j}(a)=1$ for any $j\in\{2,3,4,5\}$, 
and when $G$ is of type $\mathbf E_7$, 
$$
S_1=SL_2(\tz R)\times Spin(4,5).
$$
The group $\tilde{U}=U/\ker(\sigma)$ is a Heisenberg group and 
$S_1$ is a subgroup of the group of automorphisms
of $\tilde{U}$ which fix $\mathcal Z(\tilde{U})$
pointwise.
This means that $S_1$ acts through
a subgroup of $Sp(\tilde{U}/\mathcal Z(\tilde{U}))$.
Therefore the existence of the extension of $\sigma$ to $R_1$ is 
immediate once we show that $S_1$ acts through a subgroup of the 
metaplectic cover $Mp(\tilde{U}/\mathcal Z(\tilde{U}))$
of $Sp(\tilde{U}/\mathcal Z(\tilde{U}))$.

When $G$ is of type $\mathbf E_6$ 
the action of 
$\tz R^\times\ltimes Spin(3,4)$ 
leaves
$\germ X$ and $\germ Y$ 
invariant. In fact if $GL(\germ X)$ denotes the Levi factor
of the Siegel parabolic of $Sp(\tilde{U}/\mathcal Z(\tilde{U}))$,
then $\tz R^\times\ltimes Spin(3,4)$ is a subgroup of  
the component group of $GL(\germ X)$.
The formulas in \cite[\S 1]{howerank}
for the oscillator
representation on the Siegel parabolic
imply that 
$S_1$ acts through a subgroup of $Mp(\tilde{U}/\mathcal Z(\tilde{U}))$.

When $G$ is of type $\mathbf E_7$ the situation is only slightly more complicated.
Recall that 
$$
S_1=SL_2(\tz R)
\times Spin(4,5).
$$ 
Similar to when $G$ is of type $\mathbf E_6$, 
$Spin(4,5)$ leaves the vector spaces $\germ X$ and $\germ Y$ 
invariant. 
Therefore $Spin(4,5)$ acts through a subgroup of the component group of
$GL(\germ X)$. Consequently, 
the group homomorphism 
$$
Spin(4,5)\mapsto Sp(\tilde{U}/\mathcal Z
(\tilde{U}))
$$
breaks into the composition of two group homomorphisms
\begin{equation}
\label{embed1}
Spin(4,5)\mapsto Mp(\tilde{U}/\mathcal Z(\tilde{U}))
\mapsto Sp(\tilde{U}/\mathcal Z
(\tilde{U}))
\end{equation}

\noindent where the second map is the natural projection from the metaplectic group
to the symplectic group.

The action of $SL_2(\tz R)$ does not preserve the 
polarization $\germ X\oplus\germ Y$ of $\tilde{U}$. 
However, one can choose a different polarization
which is preserved by the action of $SL_2(\tz R)$. This polarization
can be described as follows. 
Let $\Omega$ be the set of roots $\alpha\in\Delta^+$ such that 
$\germ g_\alpha\subset\germ X\oplus\germ Y$. Then $\Omega$
can be partitioned into a disjoint union \vspace{-3mm}
$$
\Omega=\bigsqcup_{t=1}^8\Omega_t
$$\vspace{-4mm}\\
where
for any $t$, $\Omega_t$ has four elements and moreover,
elements of each $\Omega_t$ can be ordered
such that
we have\vspace{-3mm}
\begin{equation}
\label{omegatsorted}
\Omega_t=\{\alpha^{(1)},\alpha^{(2)},
\alpha^{(3)},\alpha^{(4)}\},
\end{equation}
\vspace{-5mm}\\
where the following two conditions hold:
\begin{itemize}
\item[a.] For some $j\in\{1,2\}$ depending on $t$,
we have 
$$
\alpha^{(2)}=\alpha^{(1)}+\alpha_7,
\alpha^{(3)}=\beta_j-\alpha^{(2)}\textrm{ and } 
\alpha^{(4)}=\alpha^{(3)}+\alpha_7.$$
\item[b.] $\germ g_{\alpha^{(1)}}\oplus
\germ g_{\alpha^{(3)}}\subseteq \germ X$ and
 $\germ g_{\alpha^{(2)}}\oplus
\germ g_{\alpha^{(4)}}\subseteq \germ Y$.
\end{itemize}

A polarization preserved by $SL_2(\tz R)$ has the form
$\germ X_1\oplus\germ Y_1$,
where 
for every $1\leq t \leq 8$, if $\Omega_t$ is 
sorted as in (\ref{omegatsorted}) then
$\germ X_1$ 
contains the direct sum $\germ g_{\alpha^{(1)}}\oplus\germ g_{\alpha^{(2)}}$
and 
$\germ Y_1$ contains the direct sum 
$\germ g_{\alpha^{(3)}}\oplus\germ g_{\alpha^{(4)}}$.
For the reader's convenience, 
we give one such  
polarization explicitly in section \ref{Appendix}.

Since $SL_2(\tz R)$ preserves a polarization,
we can argue as above and see that the group homomorphism
\vspace{-2mm}
$$
SL_2(\tz R)\mapsto Sp(\tilde{U}/\mathcal Z(\tilde{U}))
$$ \vspace{-5mm}\\
breaks into the composition of
two group homomorphisms 
\begin{equation}
\label{embed2}
SL_2(\tz R)\mapsto Mp(\tilde{U}/\mathcal Z(\tilde{U}))
\mapsto Sp(\tilde{U}/\mathcal Z(\tilde{U})).
\end{equation}
Next we show that the group homomorphism
\vspace{-2mm}
$$
SL_2(\tz R)\times Spin(4,5)\mapsto Sp(\tilde{U}/\mathcal Z(\tilde{U}))
$$\vspace{-5mm}\\
breaks into the composition of two group homomorphisms
\vspace{-2mm}
$$
SL_2(\tz R)\times Spin(4,5)\mapsto Mp(\tilde{U}/\mathcal Z(\tilde{U}))
\mapsto Sp(\tilde{U}/\mathcal Z(\tilde{U})).
$$ \vspace{-5mm}\\
Let $\phi_1$ and $\phi_2$ be the maps from $Spin(4,5)$ 
and $SL_2(\tz R)$ into $Mp(\tilde{U}/\mathcal Z(\tilde{U}))$ 
given in (\ref{embed1}) and (\ref{embed2}) respectively. Consider the 
map \vspace{-2mm}
$$
\Phi:Spin(4,5)\times SL_2(\tz R)\mapsto 
Mp(\tilde{U}/\mathcal Z(\tilde{U}))
$$ \vspace{-4mm}\\
given by 
$\Phi(a\times b)=\phi_1(a)\phi_2(b)$.
We show that the map $\Phi$ is the appropriate group
homomorphism from $S_1$ to $Mp(\tilde{U}/\mathcal Z(\tilde{U}))$.
Continuity of $\Phi$ is obvious.
To show that $\Phi$ is a group homomorphism, it suffices
to show that $\phi_1(a)$ and $\phi_2(b)$ 
commute. But 
the images of $\phi_1(a)$ and $\phi_2(b)$ inside 
$Sp(\tilde{U}/\mathcal Z(\tilde{U}))$ commute with each other, and 
since $Spin(4,5)\times SL_2(\tz R)$ is connected,
the commutator of $\phi_1(a)$ and $\phi_2(b)$ should be a constant
function. Checking for when $a$ and $b$ are the identity elements 
implies that this commutator is equal to the identity element
of $Mp(\tilde{U}/\mathcal Z(\tilde{U}))$. 

The uniqueness of $\hat{\sigma}$ follows 
from the fact that the group $SL_2(\tz R)\times Spin(4,5)$ is 
perfect. 

\fff

\end{proof}

\section{A class of representations of rank two}
\label{degenerate_principal_series}

For any irreducible representation $\sigma$ of a 
nilpotent simply connected Lie group, 
let
$\mathcal O_\sigma$ be the coadjoint orbit associated to $\sigma$ 
(see \cite{kirillov},\cite{kirillovbook}). 
Recall the following elementary results from Kirillov's orbit method.

\begin{Proposition}
\label{kirillovneeds}
Let $N^1\subset N^2$ be nilpotent simply connected
Lie groups and assume $N^1$ is a Lie subgroup of $N^2$
of codimension $n$.

\begin{itemize}
\item[a.] If $\sigma^1$ is an irreducible 
unitary representation of $N^1$, then
the support of $\ind_{N^1}^{N^2}\sigma^1$ lies inside
irreducible representations $\sigma$ of $N^2$
for which \vspace{-2mm}
$$
\dim(\mathcal O_\sigma)\leq 2n+\dim(\mathcal O_{\sigma^1}).
$$\vspace{-7mm}
\item[b.] If $\sigma$ is an irreducible unitary representation of
$N^2$ then the support of $\res_{N^1}^{N^2}\sigma$ lies inside
irreducible representations $\sigma^1$ of $N^1$ 
for which \vspace{-2mm}
$$
\dim(\mathcal O_{\sigma^1})\leq \dim(\mathcal O_\sigma).
$$\vspace{-7mm}

\end{itemize}
\end{Proposition}

\begin{Proposition}
\label{rankofunitary}
The
unitary principal series representations 
$\pi_\chi$ of $G$ 
are of rank two.
\end{Proposition}

\begin{proof}
Recall that $N_B=[B,B]$.
From Bruhat decomposition it follows that in the double coset space
$N_B\backslash G/\overline{P}$, $N_B\overline P$ has full measure. Therefore  
by Mackey's subgroup theorem 
$$
\res_{N_B}^G\pi_\chi=\res_{N_B}^G\ind_{\overline P}^G\chi
=\ind_{N_B\cap\overline P}^{N_B}\chi.
$$
Let $U_1=N_B\cap\overline P$. Since $U_1\subset [{\overline P},{\overline P}]$,
the restriction of $\chi$ 
to $U_1$ is the trivial representation. Therefore 
$$
\res_{N_B}^G\pi_\chi=\ind_{U_1}^{N_B}1.$$

Recall that $\beta_1$ is the highest root of $\germ g$.
Let $G_{\beta_1}$ be the one-parameter unipotent subgroup of $G$
which corresponds to $\germ g_{\beta_1}$. One can see that
$U_2=U_1G_{\beta_1}$ is actually a Lie subgroup 
of $N_B$ and $U_2\approx U_1\times G_{\beta_1}$.
If $e$ denotes the identity element of $G$, then 
$$
\ind_{U_1}^{U_2}1=\ind_{U_1}^{U_1}\hat{\otimes}\ind_{\{e\}}^{G_{\beta_1}}1=
1\hat{\otimes}L^2(G_{\beta_1},dg_{\beta_1})
$$
where $L^2(G_{\beta_1},dg_{\beta_1})$ is the left regular representation
of $G_{\beta_1}$.
Consequently, $\ind_{U_1}^{U_2}1$
is actually a direct integral of one-dimensional representations of
$U_2$. Any one-dimensional representation is associated
to a coadjoint orbit of dimension zero. 
Proposition \ref{kirillovneeds}
and 
the relation
$$\ind_{U_1}^{N_B}1=\ind_{U_2}^{N_B}\ind_{U_1}^{U_2}1$$ 
imply that 
$\ind_{U_1}^{N_B}1$ (respectively
$
\res_{N_\gg}^G\pi_\chi=\res_{N_\gg}^{N_B}\ind_{U_1}^{N_B}1
$)
is supported on irreducible 
representations of $N_B$ 
(respectively $N_\gg$)
whose coadjoint orbits have dimensions
at most twice the codimension of $U_2$ in $N_B$.
A simple calculation shows that the codimension of
$U_2$ in $N_B$ is equal to 15 (respectively 26) when $G$ is of type $\mathbf E_6$ (respectively
when $G$ is of type $\mathbf E_7$).
However, by \cite[Cor. 4.2.3]{salduke}, 
the dimension of the 
coadjoint orbit associated to
a rankable representation of $N_\gg$ of rank 3 is 32 when $G$ is of type $\mathbf E_6$
and 56 when $G$ is of type $\mathbf E_7$. But $2\times 15<32$ and $2\times 26<56$(!), which 
imply that $\pi_\chi$ should be of rank at most two. 
A comparison of $K$-types and the uniqueness of the minimal 
representation \cite[Prop 14]{salmanus} justify that $\pi_\chi$ 
is not of rank one.

\fff
\end{proof}

Proposition \ref{rankofunitary} can actually be extended to 
include the complementary series representations $\pi_s$ 
and the representation
$\pi^\circ$ which appear
when $G$ is of type $\mathbf E_7$. 
Note that the proof of Proposition \ref{rankofunitary}
is not valid anymore.

\begin{Proposition}
\label{complementary_series_theorem}
$\pi^\circ$ and $\pi_s$ (for $0\!\leq\! s\!<\! 1$) 
are of rank two.

\end{Proposition}

\begin{proof}
We use a construction of these representations
given in \cite{barchinietal}. 
First consider a complementary series representation $\pi_s$. 
Using \cite[Corollary 8.7]{barchinietal} one can describe
the $\overline P$-action of $\pi_s$
on $L^2(\germ n,\nabla(X)^{-s}d_sX)$,
where $\nabla$ is a cubic $L$-semi-invariant polynomial
and $d_sX$ is a normalization of the Lebesgue measure.
The Fourier transform gives an isometry
$$
\widehat{\ \ \ }:L^2(\germ n,d_sX)\mapsto  
L^2(\overline{\germ n},d_s\overline X)$$
where $d_s\overline X$ is an appropriate normalization
of the Lebesgue measure on $\overline{\germ n}$. 
The isometry
$f\mapsto \widehat{f\nabla^{-{s\over 2}}}$
from $L^2(\germ n,\nabla(X)^{-s}d_sX)$
to 
$
L^2(\overline{\germ n},d_s\overline X)
$
is an intertwining operator
between the action of $\overline P$ on 
$L^2(\germ n,\nabla(X)^{-s}d_sX)$ given in 
\cite[(8.5)]{barchinietal}
and its action on $L^2(\overline{\germ n},d_s\overline X)$
given as follows:
\vspace{-3mm}
\begin{eqnarray}
\label{eqnparabeqn1}
{e^{\overline Y}}\cdot {f(\overline X)}&=&f(\overline X+\overline Y)
\qquad\textrm{for every }
{\overline X},{\overline Y}\in {\overline{\germ n}}
\end{eqnarray}
\vspace{-7mm}
\begin{eqnarray}
\label{eqnparabeqn2}
{l\cdot f}({\overline X})&=&\delta_P(l)^{1\over 2}
f(l^{-1}\cdot {\overline X})
\qquad\textrm{for every }
{l\in L}, {\overline X}\in{\overline{\germ n}}.
\end{eqnarray}
\vspace{-5mm}\\
Note that these formulas are independent of $s$.
From the description of the unitary principal series
in the ``non-compact'' picture in
\cite[Chapter VII, (7.3b)]{knapp}
it follows that for any unitary character $\chi$,
${\pi_\chi}_{|{\overline P}}$ and ${\pi_s}_{|{\overline P}}$
are isomorphic. It follows that 
\vspace{-2mm}
$$
\res_P^G\pi_{\chi}=\res_P^G \pi_s.
$$
\vspace{-5mm}\\
Since $N_\gg\subset P$, the rank of $\pi_s$
is equal to the rank of $\pi_\chi$.

Next we prove that $\pi^\circ$ is of rank two. Our proof actually 
shows something slightly stronger: we prove
that the 
(spectrum of the) restriction
of $\pi^\circ$ to $U\cap N$ is multiplicity-free. 
Note that this is false for any representation of $G$ of rank
three (as it already fails for a rankable representation of 
$N_\gg$ of rank three)
but is not always true for every representation of $G$
of rank two. In fact it follows from the results of this
paper
that the only representations
of $G$ with this property are the trivial representation, the minimal
representation, and $\pi^\circ$.

We can 
conjugate $N\cap U$ by the longest Weyl element and obtain 
its ``opposite'' group, $\overline N\cap \overline U$. It
suffices to prove that the $\overline N\cap \overline U$-spectrum
of $\pi^\circ$ is multiplicity-free. We use \cite[Theorem 8.11]{barchinietal}. 

The Lie algebra of 
$N\cap U$ is equal to $\germ Y\oplus \mathcal Z(\germ u)$, 
with $\germ Y$ as in (\ref{descripuxyz}),
and
we have 
$$
\germ n=\germ Y\oplus\mathcal Z(\germ u)\oplus \germ g_{\alpha_7}.
$$
If we think of $\germ n$ as the Jordan
algebra $\mathrm{Herm}(3,\tz O_{\mathrm{split}})$, then from results in
\cite[\S 2]{barchinietal} 
it follows that 
$\germ Y$, $\mathcal Z(\germ u)$
and $\germ g_{\alpha_7}$ 
correspond to matrices of the form 
\begin{equation}
\label{2222}
\left[\begin{array}{ccc}
0&v&w\\
v^*&0&0\\
w^*&0&0
\end{array}
\right]
\qquad,\qquad 
\left[\begin{array}{ccc}
t_1&u&0\\
u^*&t_2&0\\
0&0&0
\end{array}
\right]
\qquad
\textrm{ and }
\qquad
\left[
\begin{array}{ccc}
0&0&0\\
0&0&0\\
0&0&t_3
\end{array}
\right]
\end{equation}
respectively.
By \cite[Theorem 8.11]{barchinietal}, 
the 
restriction of $\pi^\circ$ to $\overline P$ can
be realized on $L^2(\mathcal O_2,d\nu_2)$, where $d\nu_2$
is the $L$-quasi-invariant measure of $\mathcal O_2$.
If $\nabla$ is an $L$-semi-invariant
cubic polynomial on $\germ n$, then
the closure of $\mathcal O_2$ in $\germ n$ 
is equal to the
vanishing set of $\nabla$ in $\germ n$.

Write $\mathrm{Herm}(3,\tz O_{\mathrm{split}})$ as 
\begin{equation}
\label{3333}
\mathrm{Herm}(3,\tz O_{\mathrm{split}})=
\mathrm{Herm}(2,\tz O_{\mathrm{split}})
\oplus\tz O_{\mathrm{split}}\oplus\tz O_{\mathrm{split}}\oplus\tz R
\end{equation}
in a fashion compatible with the matrices shown in (\ref{2222}).
More precisely, the summands
$\mathrm{Herm}(2,\tz O_{\mathrm{split}})$,
$\tz O_{\mathrm{split}}\oplus\tz O_{\mathrm{split}}$, and
$\tz R$ 
and correspond to $\mathcal Z(\germ u),\germ Y,$ and $\germ g_{\alpha_7}$
respectively.
The decomposition (\ref{3333})
allows us
to represent points of the left hand side 
of (\ref{3333}) 
by quadruples
$\mathbf u=(u_1,u_2,u_3,u_4)$, where: 
\vspace{-2mm}
$$
u_1\in\mathrm{Herm}(2,\tz O_{\mathrm{split}}),
u_2,u_3\in \tz O_{\mathrm{split}},\textrm{ and } u_4\in\tz R.
$$\vspace{-5mm}\\
From the description of the action
of $\overline N$ on $L^2(\mathcal O_2,d\nu_2)$ 
in \cite[Lemma 8.10]{barchinietal}
it follows that $\overline N\cap \overline U$
acts on a space of functions on $\mathcal O_2$ by pointwise multiplication
by characters, and two distinct points $\mathbf u$ and $\mathbf u'$ are
separated by these characters unless we have $u_i=u_i'$ for any
$i\in\{1,2,3\}$, but $u_4\neq u_4'$. 
If we can show that for a set $\mathcal S\subseteq \mathcal O_2$ of 
full measure, any two distinct 
points differ in at least one of the
first three coordinates, then it follows that 
the action of $\overline N\cap \overline U$ separates points of
$\mathcal S$, and therefore the (spectrum of the) 
action of $\overline N\cap\overline U$
on $L^2(\mathcal O_2,d\nu_2)$
is multiplicity-free. Our next aim is to prove
the existence of the set $\mathcal S$.

For an element $\mathbf u=(u_1,u_2,u_3,u_4)$, we have
$$
\nabla(\mathbf u)=u_4\nabla_1(u_1)+F(u_1,u_2,u_3)
$$ 
where $\nabla_1$ is the ``determinant'' of the Jordan algebra
$\mathrm{Herm}(2,\tz O_{\mathrm{split}})$ and $F(u_1,u_2,u_3)$ 
is a cubic polynomial in 26 variables 
obtained by real coordinates of $u_1,u_2,u_3$
(see
\cite[\S 5, (5.11)]{tonnyspringer}). 
Obviously, if $\mathbf u\in\mathcal O_2$ is
such that $\nabla_1(u_1)\neq 0$, then the equation $\nabla(\mathbf u)=0$
uniquely determines $u_4$ in terms of $u_1,u_2,u_3$. Therefore 
we can choose $\mathcal S$ to be the set of all 
$\mathbf u\in \germ n$ 
for which
$\nabla(\mathbf u)=0$ but $\nabla_1(u_1)\neq 0$. It remains to
show that this set has full measure in $\mathcal O_2$. To this end,
we show that the complement of $\mathcal S$ in $\mathcal O_2$ is a
submanifold of $\mathcal O_2$ of positive codimension. In fact we
can work with the complexifications. Let 
$\overline{\mathcal O_2}$ denote the Zariski closure of 
$\mathcal O_2$
in 
$\mathrm{Herm}(3,\tz O_{\mathrm{split}})\otimes\tz C$.
Recall that the closure of
$\mathcal O_2$ in $\germ n$ is equal to the
set of $\tz R$-rational points of the vanishing set
of $\nabla$ in
$\mathrm{Herm}(\tz O_{\mathrm{split}})\otimes \tz C$.
The latter vanishing set is the closure of 
an orbit of the action
of the complexification of $L$ on a 27-dimensional 
complex affine space, and therefore it is an 
irreducible algebraic variety. 
Moreover, there exist elements $\mathbf u$ 
of $\mathcal O_2$ for which $\nabla_1(u_1)\neq 0$. 
Therefore the set of all $\mathbf u$ in $\overline{\mathcal O_2}$
for which $\nabla_1(u_1)=0$ is an algebraic set 
of positive codimension in $\overline{\mathcal O_2}$.
Since the complex dimension of $\overline{\mathcal O_2}$ is the 
same as the real dimension of $\mathcal O_2$,
the complement of $\mathcal S$ in $\mathcal O_2$ is a
submanifold of positive codimension.

\fff

\end{proof}

\section{Proof of Theorem \ref{main2}a}

In this section we study the correspondence (\ref{correspondence})
for representations $\pi_\chi$ when $G$ is of type $\mathbf E_7$. 
Our main tool is standard Mackey theory.
Recall that $\overline P$ is the parabolic opposite to $P$, 
and $\pi_\chi=\ind_{\overline P}^G\chi$ where 
$\chi$ is a unitary character of $\overline P$. 
The parabolic $R$ can be expressed as
\begin{equation}
\label{rreminder}
R=(\tz R^\times\ltimes (SL_2(\tz R)\times Spin(5,5)))\ltimes U
\end{equation}
where $\tz R^\times$ is an appropriate subgroup of $A$.
(Conjugation by $-1\in\{\pm 1\}\subset\tz R^\times$
induces a nontrivial automorphism of $Spin(5,5)$.)
Let $\BBB$ denote the Borel subgroup of $SL_2(\tz R)$
which contains $A\cap SL_2(\tz R)$ and the unipotent subgroup corresponding to 
$-\alpha_7$. Let 
$
\NNN=[\BBB,\BBB].
$
Observe that the vector space $\germ X$
is in fact a commutative Lie subalgebra of $\germ g$ which lies
inside the Lie algebra of $N_B$. 
Therefore $\germ X$ is the Lie algebra of a Lie subgroup of
$G$. We abuse our notation slightly to denote this Lie subgroup by
$\germ X$ as well.
 
Bearing in mind that $\tz R^\times$ represents a specific 
subgroup of $A$ as in (\ref{rreminder}), 
we consider the following subgroups of $R$:
\vspace{-2mm}
\begin{eqnarray}
R_2=(\tz R^\times \ltimes(\BBB\times Spin(5,5)))\ltimes \germ X\\
R_3=(\tz R^\times \ltimes(\BBB\times Spin(5,5)))\ltimes U.
\end{eqnarray}
\vspace{-5mm}\\
Note that indeed $R_2=R\cap \overline P$.

By Bruhat decomposition, $R\overline P$ is an open double coset
of full measure in $R\backslash G/{\overline P}$. Therefore Mackey's  
subgroup theorem implies that
\vspace{-2mm}
\begin{equation}
\label{hasan}
\res_R^G{\pi_\chi}=\ind_{R_2}^R\chi.
\end{equation}
\vspace{-5mm}\\
Next observe that the isomorphism
\vspace{-2mm}
$$
R_3/((\NNN\times Spin(5,5))\ltimes U)
\approx 
R_2/((\NNN\times Spin(5,5)\ltimes \germ X)
$$
\vspace{-4mm}\\
implies that $\chi$ extends to a character $\hat{\chi}$
of $R_3$.
By the projection formula
\vspace{-1mm}
\begin{equation}
\label{gholumi}
\ind_{R_2}^R\chi=\ind_{R_3}^R\ind_{R_2}^{R_3}\chi=
\ind_{R_3}^R(\hat{\chi}\otimes\ind_{R_2}^{R_3}1).
\end{equation}

Let $\eta=\ind_{R_2}^{R_3}1$. Let $\sigma$ be the
big representation of $U$ introduced in the statement
of Proposition \ref{extensionispossible}
and let $\hat{\sigma}$ be the 
extension of $\sigma$ to $R_1$.
We will now prove the following lemma.
\begin{Lemma}
\label{lemmafore7} 
There exists a unitary character $\hat{\kappa}$
of
$R_3$ such that
\vspace{-2mm}
\begin{equation}
\label{keylemmainduction}
\eta=\hat{\kappa}\otimes(\res_{R_3}^R\ind_{R_1}^R\hat{\sigma}).
\end{equation}
\end{Lemma}
\begin{proof}
To prove Lemma \ref{lemmafore7} note that
$UR_2=R_3$ and therefore
by Mackey's subgroup theorem\vspace{-3mm}
$$
\res_U^{R_3}\eta=\res_U^{R_3}\ind_{R_2}^{R_3}1=\ind_{\germ X}^U1.
$$
\vspace{-4mm}\\
The right hand side, which is equal to the $U$-spectrum of 
$\eta$, 
is a multiplicity-free direct integral of
big representations of $U$. The fact that the action of 
$\tz R^\times\ltimes Spin(5,5)\subset R_3$ on the $U$-spectrum of $\eta$ is
transitive implies that $\eta$ is an 
irreducible representation of $R_3$. Next 
we apply standard Mackey theory to $\eta$.
The stabilizer
of $\sigma$ in $R_3$ is $R_3\cap R_1$, and therefore by Mackey theory
we can write
\vspace{-2mm}
\begin{equation}
\label{labadoo}
\eta=\ind_{R_1\cap R_3}^{R_3}(\kappa\otimes \hat{\sigma})
\end{equation}
\vspace{-4mm}\\
where $\kappa$ is an irreducible unitary representation of 
$\BBB\times Spin(4,5)$ which is extended (trivially on $U$) 
to $R_1\cap R_3$.
We can use Mackey's subgroup theorem again and 
see that
the $U$-spectrum of the right hand side of (\ref{labadoo})
is multiplicity-free only if
$\kappa$ is a one-dimensional representation. Therefore 
$\kappa$ should be a unitary character of 
$\BBB\times Spin(4,5)$ and hence it is trivial on $\NNN\times Spin(4,5)$.
Obviously $\kappa$ extends to a unitary character $\hat{\kappa}$
of 
$R_3$ which factors through a character of $R_3/((N_{SL_2}\times Spin(5,5))\ltimes U)$. 
From (\ref{labadoo}), the projection formula, Mackey's subgroup theorem, and $R_3R_1=R$
we have
$$
\eta=\hat{\kappa}\otimes\ind_{R_1\cap R_3}^{R_3}\hat{\sigma}=
\hat{\kappa}\otimes \res_{R_3}^R\ind_{R_1}^R\hat{\sigma}.
$$

\fff

\end{proof}

Let $\hat{\eta}=\ind_{R_1}^R\hat{\sigma}$. From  
(\ref{keylemmainduction}) it follows that 
\begin{equation}
\label{wwwghol}
\hat{\kappa}^{-1}\otimes\eta=\res_{R_3}^R\hat{\eta}.
\end{equation}
Next we continue with the right hand side of 
(\ref{gholumi}). Using (\ref{wwwghol}) and the projection formula
we have 
\begin{equation}
\label{7.7.7}
\ind_{R_3}^R(\hat{\chi}\otimes \eta)=\ind_{R_3}^R(\hat{\chi}\otimes\hat{\kappa}
\otimes
\hat{\kappa}^{-1}\otimes\eta)
=\hat{\eta}\otimes\ind_{R_3}^R(\hat{\kappa}\otimes\hat{\chi})
=\hat{\eta}\otimes\hat{\zeta}
\end{equation}
where
$\hat{\zeta}=\ind_{R_3}^R (\hat{\kappa}\otimes\hat{\chi})$. 
By Mackey's subgroup theorem
$$
\res_{SL_2(\tz R)}^R\hat{\zeta}=\res_{SL_2(\tz R)}^R\ind_{R_3}^R
(\hat{\kappa}\otimes\hat{\chi})
=\ind_{\BBB}^{SL_2(\tz R)}(\hat{\kappa}\otimes\hat{\chi}).
$$
Let $\zeta=\ind_{\BBB}^{SL_2(\tz R)}(\hat{\kappa}\otimes\hat{\chi})$.
Clearly $\zeta$ is a unitary 
principal series 
representation of $SL_2(\tz R)$. 
Since the only one-dimensional representation of $Spin(4,5)$ 
is the trivial representation, it follows immediately
that 
\begin{equation}
\label{goozgooz}
\res_{SL_2(\tz R)\times Spin(4,5)}^R\hat{\zeta}=
\zeta\hat{\otimes}1.
\end{equation}
By Mackey's subgroup theorem
we see that the representation
$$
\res_{U}^R \hat{\zeta}=\res_{U}^R\ind_{R_3}^R(\hat{\kappa}\otimes\hat{\chi})
$$
is a direct integral of representations of the form $\ind_U^U1$, i.e.,
$\hat{\zeta}$ acts trivially when restricted to $U$, or in other words
$\hat{\zeta}$ factors through a representation of $R/U$.

Let and $\tilde{\zeta}=\res_{R_1}^R\hat{\zeta}$. 
Then by (\ref{hasan}), (\ref{gholumi}), (\ref{7.7.7}), and the projection
formula we have
$$
\res_R^G\pi_\chi=\ind_{R_3}^R(\hat{\chi}\otimes\eta)=
\hat{\eta}\otimes \hat{\zeta}=(\ind_{R_1}^R\hat{\sigma})
\otimes \hat{\zeta}
=\ind_{R_1}^R(\hat{\sigma}\otimes \tilde{\zeta}).
$$
Since $\tilde{\zeta}$ comes from a representation
of $R_1/U$, 
it follows from standard Mackey theory and (\ref{goozgooz})
that $\Psi(\pi_\chi)=\zeta\hat{\otimes}1$.

\section{The image of $\pi^\circ$ under $\Psi$}
From the
$K$-type structure of the representation
$\pi^\circ$ it follows that its Gelfand-Kirillov dimension is
strictly larger than the minimal representation but
strictly smaller than a generic degenerate principal series representation
(see \cite{sahiinv}). 
This suggests that in the correspondence (\ref{correspondence}),
the image of $\pi^\circ$ should be a representation of 
$SL_2(\tz R)\times Spin(4,5)$ which is ``smaller than'' the image of $\pi_\chi$.
It turns out that the only 
possibility is the trivial representation of $SL_2(\tz R)\times Spin(4,5)$.
A rigorous proof of this statement can be given using the
property of $\pi^\circ$ which was shown 
in the proof of Proposition  
\ref{complementary_series_theorem}:
since the $N\cap U$-spectrum of $\pi^\circ$ is multiplicity-free,
Mackey's subgroup theorem applied to
$$
\res_{N\cap U}^R\ind_{R_1}^R(\Psi(\pi^\circ)\otimes\hat{\sigma})
$$
implies 
that $\Psi(\pi^\circ)$ should be one-dimensional. But 
the only one-dimensional representation of $SL_2(\tz R)\times Spin(4,5)$ is
the trivial representation. In other words, $\Psi(\pi^\circ)$ is the
trivial representation.

\section{Proof of Theorem \ref{main2}b}

In this section we study the correspondence of (\ref{correspondence})
for representations $\pi_\chi$ when $G$ is of type $\mathbf E_6$. 
The argument is very similar to the case when $G$ is of type $\mathbf E_7$.

Recall that $\overline{P}$ is the parabolic subgroup opposite to $P$
and $\chi$ is a unitary multiplicative character of $\overline P$ 
such that $\pi_\chi=\ind_{\overline P}^G\chi$. 

The vector space $\germ X$ is a Lie algebra, and corresponds to a Lie subgroup
of $N_B$. We abuse our notation to
let $\germ X$ denote this Lie subgroup of $N_B$ as well.
Let $R_2=R\cap \overline P$. Then
$$
R_2=((\tz R^\times\times\tz R^\times)\ltimes Spin(4,4))\ltimes \germ X
$$ 
Since the restriction of $\chi$ to $Spin(4,4)\ltimes \germ X$ 
is trivial, $\chi$ extends to a character $\hat{\chi}$
of $R$. By the projection formula
$$
\res_R^G\ind_{\overline P}^G\chi=\ind_{R\cap \overline P}^R\chi
=\hat{\chi}\otimes \ind_{R\cap \overline P}^R1.
$$

\begin{Lemma}
\label{lemmafore6}
Fix an extension $\hat{\sigma}$ of $\sigma$ to $R_1$. 
(See the statement of Proposition
\ref{extensionispossible}.) Then
for some unitary character $\tilde{\chi}$ of $R_1$ whose restriction to
$U$ is trivial, we have
\begin{equation}
\label{labudi}
\ind_{R\cap \overline P}^R1=
\ind_{R_1}^R(\tilde{\chi}\otimes \hat{\sigma}).
\end{equation}
\end{Lemma}
 
\begin{proof}
By Mackey's subgroup theorem, we have
$$
\res_U^R\ind_{R\cap \overline P}^R 1=\ind_{\germ X}^U1
$$
which implies that the restriction of the left hand side of 
(\ref{labudi}) to $U$ is a multiplicity-free direct integral of
big representations of $U$. Transitivity of the action of 
the Levi subgroup of $R$ on this $U$-spectrum
implies that the left hand side of (\ref{labudi})
is an irreducible representation of $R$.
It follows that if we use standard Mackey theory to
write the left hand side of (\ref{labudi}) as 
$\ind_{R_1}^R(\tau\otimes\hat{\sigma})$, 
then $\tau$ has to be a one-dimensional
representation of $R_1/U$. 
Therefore $\tau=\tilde{\chi}$, for some unitary character
$\tilde{\chi}$.

\fff

\end{proof}

Using Lemma \ref{lemmafore6}, we have
$$
\res_R^G\pi_\chi=\hat{\chi}\otimes \ind_{R\cap \overline P}^R1=
\hat{\chi}\otimes \ind_{R_1}^R(\tilde{\chi}\otimes \hat{\sigma})=
\ind_{R_1}^R(\hat{\chi}\otimes \tilde{\chi}\otimes\hat{\sigma})
$$
which implies that $\Psi(\pi_\chi)=\tilde\chi\otimes\hat{\chi}$.

\section{Tables}
\label{Appendix}

In this section we prove Lemma \ref{basislemma} 
and describe explicitly the 
polarization $\germ X_1\oplus\germ Y_1$
used in the proof of Proposition 
\ref{extensionispossible}.
Recall the labelling of simple roots shown in section \ref{notation}. 
For any root $\alpha$,
we can write 
$$
\alpha=c_1(\alpha)\alpha_1+\cdots+c_n(\alpha)\alpha_n
$$
for integers $c_i(\alpha)$. We represent a root by putting these integers
in their corresponding locations on the nodes of the Dynkin diagram.
When describing a set of roots, a ``*'' in a location of the 
Dynkin diagram means that the corresponding coefficient $c_i$ 
can assume any possible value, given the fixed coefficients, to
make the entire labelling represent a root.

Parts a,b,c and e of Lemma \ref{basislemma} follow easily from the
given tables. Part d requires a simple calculation of weights
of the standard basis of $\tz R^{n-2,n-2}$ under the action of 
the Cartan subgroup of $Spin(n-2,n-2)$.

Let $G$ be of type $\mathbf E_6$. 
Tables below suffice for the proof of Lemma
\ref{basislemma}.\\

\begin{center}
\begin{tabular}{|c|c|c|c|c|}
\hline
$\germ W$ & $\germ W^*$ &  $\germ X$
 & $\germ Y$ &$\mathcal Z(\germ u)$ \\
\hline

\jad{0****}{1}

&

\jad{1****}{1}

&\jad{0***1}{*} &\jad{1***0}{*} &
\jad{1***1}{*}\\

\hline
\end{tabular}
\vspace{5mm}

\begin{tabular}{|c|c|c|c|}
\hline 
$e_1$  &  $e_2$  &  $e_3$ &  $e_4$ \\

\hline
\jad{12321}{2} & \jad{12321}{1} & \jad{12221}{1} & \jad{12211}{1}\\

\hline\hline
$e_{-1}$ & $e_{-2}$ &  
$e_{-3}$  & $e_{-4}$\\
\hline

\jad{11111}{0} & \jad{11111}{1} & \jad{11211}{1} & \jad{11221}{1}\\

\hline
\end{tabular}
\vspace{5mm}

\begin{tabular}{|c|c|c|c|c|c|}
\hline
  $f_2$  &  $f_3$ &  $f_4$ & $f_{-2}$ & $f_{-3}$ &  
$f_{-4}$  \\

\hline

 \jad{00000}{1} & \jad{00100}{1} & \jad{00110}{1} &\jad{01210}{1} &
\jad{01110}{1} & \jad{01100}{1} \\
\hline 
\end{tabular}
\vspace{4mm}

\end{center}

Now let $G$ be of type $\mathbf E_7$. 
Tables below suffice for the proof of 
Lemma \ref{basislemma}.\\

\begin{center}
\footnotesize
\begin{tabular}{|c|c|c|c|c|}
\hline

$\germ W$ & $\germ W^*$ &  $\germ X$  & $\germ Y$ &$\mathcal Z(\germ u)$\\

\hline 

\jadv{0****1}{\:\:*} & \jadv{1****1}{\:\:*} & \jadv{01****}{*} & \jadv{11****}{*}&
\jadv{12****}{\:*}\\

\hline

\end{tabular}

\vspace{10mm}

\begin{tabular}{|c|c|c|c|c|}
\hline 
$e_1$  &  $e_2$  &  $e_3$ &  $e_4$ & $e_5$\\

\hline 

\jadv{123432}{2} & \jadv{123431}{2} & \jadv{123421}{2} &
\jadv{123321}{2} &\jadv{123321}{1}\\
 
\hline 
\hline
$e_{-1}$ & $e_{-2}$ &  
$e_{-3}$  & $e_{-4}$ &$e_{-5}$\\
\jadv{122210}{1}&
\jadv{122211}{1} &\jadv{122221}{1} &\jadv{122321}{1} &
\jadv{122321}{2}\\
\hline
\end{tabular}

\vspace{5mm}

\begin{tabular}{|c|c|c|c|}

\hline

$f_2$  &  $f_3$ &  $f_4$ & $f_5$\\
\hline

\jadv{000001}{0} &\jadv{000011}{0} & \jadv{000111}{0} &
\jadv{000111}{1}   \\

\hline\hline
 $f_{-2}$ & $f_{-3}$ &  
$f_{-4}$  & $f_{-5}$ \\
\hline
\jadv{001221}{1} & \jadv{001211}{1} &
\jadv{001111}{1} & \jadv{001111}{0}\\

\hline

\end{tabular}
\vspace{5mm}
\end{center}

Finally, when $G$ is of type $\mathbf E_7$
the polarization $\germ X_1\oplus\germ Y_1$
which appears in the proof of Proposition 
\ref{extensionispossible} can be described as follows.

\begin{center}
\begin{tabular}{|c|c|}
\hline
$\germ X_1$&$\germ Y_1$\\
\hline
\begin{tabular}{c}
\jadv{*1**00}{\,\:\:*}\\
\jadv{*1**11}{\,\:\:*}
\end{tabular}
&\begin{tabular}{c}
\jadv{*1**10}{\,\:\:*}\\
\jadv{*1**21}{\,\:\:*}
\end{tabular}\\
\hline
\end{tabular}
\end{center}

\end{document}